\documentclass[12pt]{amsart}

\usepackage{amsthm, amssymb, amsmath, amsfonts}
\usepackage{enumerate, enumitem}

\usepackage[colorlinks]{hyperref}

\allowdisplaybreaks

\newtheorem{thm}{Theorem}[section]
\newtheorem{thms}{Theorem}[subsection]
\newtheorem{lems}{Lemma}[subsection]

\newtheorem{rem}[thm]{Remark}
\newtheorem{rems}[thms]{Remark}
\newtheorem*{nota}{Notation}
\newtheorem{thmss}[lems]{Theorem}

\title[Restriction theorems for $B_p(G)$]{Restriction theorems for the $p$-analog of the Fourier-Stieltjes algebra}
\author[A. Dabra]{Arvish Dabra}
\address{Arvish Dabra,\newline\indent Department of Mathematics,\newline\indent Indian Institute of Technology Delhi,\newline\indent New Delhi - 110016, India.}
\email{arvishdabra3@gmail.com}
\author[N. S. Kumar]{N. Shravan Kumar$^\ast$}
\address{N. Shravan Kumar,\newline\indent Department of Mathematics,\newline\indent Indian Institute of Technology Delhi,\newline\indent New Delhi - 110016, India.}
\email{shravankumar.nageswaran@gmail.com}

\begin{document}
	
	\keywords{locally compact groups, representations, QSL$_p$-spaces, $p$-analog of Fourier-Stieltjes algebra, restriction. \\
 $\ast$ Corresponding author: shravankumar.nageswaran@gmail.com}
	

	\maketitle
	
	\begin{abstract}
		For a locally compact group $G$ and $1 < p < \infty,$ let $B_p(G)$ denote the $p$-analog of the Fourier-Stieltjes algebra $B(G) \, (\text{or} \, B_2(G))$. Let $r: B_p(G) \to B_p(H)$ be the restriction map given by $r(u) = u|_H$ for any closed subgroup $H$ of $G.$ In this article, we prove that the restriction map $r$ is a surjective isometry for any open subgroup $H$ of $G.$ Further, we show that the range of the map $r$ is dense in $B_p(H)$ when $H$ is either a compact normal subgroup of $G$ or compact subgroup of an [SIN]$_H$-group. 
	\end{abstract}
	
	\section{Introduction}
	For a locally compact group $G,$ Eymard \cite{Eym} introduced and studied the Fourier and Fourier-Stieltjes algebra, denoted by $A(G)$ and $B(G),$ respectively. When the group $G$ is abelian, the algebras $A(G)$ and $B(G)$ are isometrically isomorphic to $L^1(\widehat{G})$ and $M(\widehat{G}),$ respectively, where $M(\widehat{G})$ is the measure algebra of the dual group $\widehat{G}.$ These Banach algebras are the central objects of study in abstract harmonic analysis. The functorial properties of these algebras were extensively studied by Eymard \cite{Eym}.
	\par
	For a closed subgroup $H$ of $G,$ consider the restriction map $r: B(G) \to B(H)$ given by $r(u) = u|_H.$ It is well known \cite{herz1973} that the map $r$ is an onto isometry for the Fourier algebra $A(G).$ However, the surjectivity of the restriction map is not true in general for the Fourier-Stieltjes algebra. For $B(G),$ it is shown in \cite[\textsection2]{mc} that $B(G)|_H = B(H)$ if $H$ is open or compact subgroup of $G.$ In \cite[Prop. 1.1]{mis}, the surjectivity of the map $r$ is established when the subgroup $H$ is the center or the identity component of $G.$ Cowling and Rodway \cite{cow} proved this result when $H$ is either a closed normal subgroup of $G$ or closed subgroup of an [SIN]$_H$-group.
	\par
	For $1 < p < \infty,$ Fig\`{a}-Talamanca \cite{figa} introduced the $p$-analog of the Fourier algebra $A(G)$ by replacing $L^2(G)$ with $L^p(G)$ when the group $G$ is abelain or compact. In 1971, Herz \cite{herz1971} generalized this algebra for any locally compact group $G.$ This algebra is denoted by $A_p(G)$ and termed as Fig\`{a}-Talamanca Herz algebra. For $p=2,$ $A_2(G)$ coincides with $A(G).$ In 2005, Runde \cite{runde} introduced and studied the $p$-analog of the Fourier-Stieltjes algebra, denoted by $B_p(G).$ Over the last two decades, even though several studies have been made for the algebra $B_p(G),$ many of its functorial properties are still not known.
	\par
	The article is organized as follows. In Section \ref{sec2}, we recall some basics of representations of a locally compact group $G$ on a Banach space $E.$ In this article, we work with a particular class of Banach spaces, namely, QSL$_p$-spaces. The related definitions and some basic properties of these Banach spaces are also recalled in this section. In Section \ref{sec4}, we prove that the restriction map $r: B_p(G) \to B_p(H)$ is a surjective isometry for any open subgroup $H$ of $G.$ Further, we show that the range of the map $r$ is dense in $B_p(H)$ when $H$ is either a compact normal subgroup of $G$ or compact subgroup of an [SIN]$_H$-group.
	
	\section{Preliminaries}\label{sec2}
	We begin this section with the definition of a representation of $G$ on a Banach space $E.$ The notations and definitions in this section are in accordance with Runde \cite{runde}.
	
	Let $\mathcal{I}(E)$ denote the set of invertible isometries on $E.$ A {\it representation} of a group $G$ on a Banach space $E$ is a pair $(\pi,E)$ where $\pi:G \to \mathcal{I}(E)$ and the map $\pi$ is a group homomorphism which is continuous with respect to the given topology on $G$ and the strong operator topology (equivalently, weak operator topology \cite[Theorem 2.8]{lee}) on $\mathcal{I}(E) \subseteq \mathcal{B}(E).$
	
	\begin{rem}
		Any representation $(\pi,E)$ of $G$ induces a representation of the group algebra $L^1(G)$ to $\mathcal{B}(E)$ given by
		$$\pi(f):= \int\limits_G f(x)\pi(x)dx \hspace{1 cm} (f \in L^1(G)),$$
		where the integral converges with respect to the strong operator topology.
	\end{rem}
	
	Two representations $(\pi,E)$ and $(\sigma,F)$ of $G$ are said to be {\it equivalent} if there is an invertible isometry $W:E \to F$ such that
		$$W \pi(x) W^{-1} = \sigma(x) \hspace{1cm} (x \in G).$$
	Throughout this article, any two representations are considered to be same if they are equivalent.
	
	As mentioned earlier, we shall be working with a particular class of Banach spaces, namely, QSL$_p$-spaces. Let $1 < p < \infty.$
	\renewcommand{\theenumi}{\roman{enumi}}
	\begin{enumerate}
		\item A Banach space is called an {\it $L_p$-space} if it is of the form $L^p(X)$ for some measure space $X.$
		\item A Banach space is called a {\it QSL$_p$-space} if it is isometrically isomorphic to a quotient of a subspace of an $L_p$-space.
	\end{enumerate}
	
	The following are some basic properties of QSL$_p$-spaces.
	\renewcommand{\theenumi}{\roman{enumi}}
	\begin{enumerate}
		\item A Banach space is a QSL$_p$-space if and only if it is a subspace of a quotient of an $L_p$-space.
		\item By definition, the class of QSL$_p$-spaces is closed under taking subspaces and quotients.
		\item The class of QSL$_p$-spaces is closed under $\ell_p$-direct sum. More precisely, if $(X_\alpha)_{\alpha \in \Lambda}$ is a family of QSL$_p$-spaces, then $\ell_p-\bigoplus\limits_{\alpha \in \Lambda} E_\alpha$ is again a QSL$_p$-space.
		\item Every QSL$_p$-space is reflexive.
	\end{enumerate}
	
	For $1 < p < \infty,$ the $p$-analog of the Fourier-Stieltjes algebra $B_p(G),$ introduced by Runde \cite{runde}, is defined as the space of all coefficient functions of isometric representations of $G$ on QSL$_{p'}$-spaces, where $1/p + 1/p' = 1.$ Let $(\pi,E)$ be a representation of $G.$ A {\it coefficient function} of $(\pi,E)$ is a function $f: G \to \mathbb{C}$ such that
	$$f(x) = \langle \pi(x)\xi,\phi \rangle \hspace{1 cm} (x \in G),$$
	where $\xi \in E$ and $\phi \in E^*.$ It is easy to verify that every coefficient function is continuous and bounded. The collection of all (equivalent classes) of representations of $G$ on a QSL$_p$-space is denoted by Rep$_p(G).$ The algebra $B_p(G)$ is defined as
	$$B_p(G):= \left\{f:G \to \mathbb{C}, f \, \text{is a coefficient function of some}\, (\pi,E) \in \text{Rep}_{p'}(G)\right\}.$$
	
	Let $(\pi_n, E_n)_{n \in \mathbb{N}} \in \text{Rep}_{p'}(G)$ be a countable family of representations. As each $E_n$ is a QSL$_{p'}$-space and since they are closed under $\ell_{p'}$-direct sum, one can consider another QSL$_{p'}$-space $E := \ell_{p'}-\bigoplus\limits_{n \in \mathbb{N}} E_n$. Now, for $\pi := \bigoplus\limits_{n \in \mathbb{N}} \pi_n,$ the pair $(\pi, E) \in \text{Rep}_{p'}(G).$ Then, from \cite[Lemma 2.3]{runde}, the function $f \in B_p(G),$ where
	$$f(x) := \sum_{n \in \mathbb{N}} \langle \pi_n(x)\xi_n,\phi_n \rangle \hspace{1cm} (x \in G),$$
	with $\xi_n \in E_n$ and $\phi_n \in {E_n}^*$ for each $n \in \mathbb{N}.$
	
	The norm on $B_p(G)$ is defined using cyclic representations. A representation $(\pi,E)$ is said to be {\it cyclic} if there is $v \in E$ such that $\pi(L^1(G))v$ is dense in $E.$ The set Cyc$_{p'}(G)$ consists of all $(\pi,E) \in \text{Rep}_{p'}(G)$ such that $(\pi,E)$ is cyclic. By \cite[Definition 4.2]{runde}, for $f \in B_p(G),$ the norm $\|\cdot\|$ is defined as the {\it infimum} over all the expressions $$\sum\limits_{n \in \mathbb{N}} \|\xi_n\| \, \|\phi_n\|,$$ where for each $n \in \mathbb{N},$ there is $(\pi_n,E_n) \in \text{Cyc}_{p'}(G)$ with $\xi_n \in E_n$ and $\phi_n \in {E_n}^*$ such that
	$$\sum_{n\in \mathbb{N}} \|\xi_n\| \, \|\phi_n\| < \infty \hspace{1cm} \text{and} \hspace{1 cm} f(x) = \sum_{n \in \mathbb{N}} \langle \pi_n(x)\xi_n,\phi_n \rangle \hspace{0.5cm} (x \in G).$$
	It is proved in \cite{runde} that $B_p(G)$ is a commutative Banach algebra that contractively contains the Fig\`{a}-Talamanca Herz algebra $A_p(G).$
	
	Throughout the article, $G$ denote a locally compact group.
	
	\section{Restriction Theorems}\label{sec4}
	Let $H$ be a closed subgroup of $G.$ In this section, we study the restriction map $r: B_p(G) \to B_p(H),$ for certain classes of closed subgroups of $G.$ Specifically, we demonstrate in Subsection \ref{sec41} that the map $r$ is a surjective isometry for each open subgroup $H$ of $G.$ Furthermore, it is proved in Subsections \ref{sec42} and \ref{sec43}, respectively, for every compact normal subgroup of a group $G$ and for any compact subgroup $H$ of an [SIN]$_H$-group that the range of the restriction map $r$ is dense in $B_p(H).$
	
	\subsection{For an open subgroup $H$ of $G$}\label{sec41}
	Let $G$ be a locally compact group and $H$ be an open subgroup of $G.$ Then, $H$ is closed and $G/H$ is a discrete space. The counting measure on $G/H$ is an invariant Haar measure.
	\par
	We first recall the $p$-induction of representations on Banach spaces, given by Jaming and Moran \cite{jaming}. The theory of $p$-induced representations along with the ideas of \cite{kaniuth} on inducing representations from an open subgroup (in case of Hilbert spaces) plays a pivotal role in proving the main result of this subsection. Let us recall the induction from an open subgroup in the setting of QSL$_{p'}$-spaces.
	\par
	Let $(\pi,E) \in \text{Rep}_{p'}(H),$ i.e., $\pi$ is a strongly continuous isometric representation of the open subgroup $H$ of $G$ in a QSL$_{p'}$-space $E.$ For the pair $(\pi,E),$ the associated Banach space $L^{p'}(G,H,\pi,E)$ is defined as
	$$L^{p'}(G,H,\pi,E) := \left\{ f:G \to E \, , \begin{array}{c} f(xh) = \pi(h)^{-1}f(x) \, \forall \, x\in G, h \in H\\\text{and} \sum\limits_{xH \in G/H} \|f(x)\|^{p'} < \infty \end{array}\right\},$$
	with the norm $$\|f\|_{p'} := \left( \sum\limits_{xH \in G/H} \|f(x)\|^{p'} \right)^{1/{p'}}.$$
	As mentioned in \cite{jaming}, the induced representation Ind$_{H \uparrow G}(\pi)$ is the representation of $G$ on $L^{p'}(G,H,\pi,E)$ and is given by 
	$$(\text{Ind}_{H \uparrow G}(\pi))(x)f(y) := (L_x(f))(y) = f(x^{-1}y).$$
	
	Motivated from \cite{kaniuth}, we consider the following Banach space:
	$$\ell^{p'}(G/H,E) := \left\{g:G/H \to E \, , \sum\limits_{v \in G/H} \|g(v)\|^{p'} < \infty\right\},$$ with the norm $$\|g\|_{p'} := \left( \sum\limits_{v \in G/H} \|g(v)\|^{p'} \right)^{1/{p'}}.$$
	
	As $G/H$ is discrete, any choice $\gamma: G/H \to G$ of a section of the $H$-cosets is a regular section, i.e., $\gamma$ is a morphism and $q \circ \gamma = I_{G/H},$ where $q$ is the canonical quotient map. For a fix choice of $\gamma,$ for $x \in G,$ $$x = \gamma(q(x)) [\gamma(q(x))^{-1} x]$$ is the unique way of writing $x$ in the form $\gamma(v)h$ with $v \in G/H$ and $h \in H.$ Then, there is a map M (depending on the choice of $\gamma$) from $\ell^{p'}(G/H,E)$ to $L^{p'}(G,H,\pi,E)$ defined by
	$$Mf(x) = Mf(\gamma(v)h) = \pi(h)^{-1}f(v) \hspace{1cm} (f \in \ell^{p'}(G/H,E)).$$
	It is easy to verify that the map $M$ is an isometric isomorphism. Further, the inverse map $M^{-1}: L^{p'}(G,H,\pi,E) \to \ell^{p'}(G/H,E)$ is given by $M^{-1}f = f \circ \gamma.$
	\par
	Now, with the help of this map $M,$ one can define a representation of $G$ on $\ell^{p'}(G/H,E).$ For $x \in G,$ let $$U_\gamma^\pi(x) := M^{-1}(\text{Ind}_{H \uparrow G}(\pi)(x))M.$$ By \cite[Definition 1.2]{runde}, it follows that $U_\gamma^\pi$ and Ind$_{H \uparrow G}$ are equivalent.
	
	\begin{rems}\label{QSP}
		Observe that $\ell^{p'}(G/H,E) = \ell_{p'}-\bigoplus\limits_{v \in G/H}E$ is a QSL$_{p'}$-space as QSL$_{p'}$-spaces are closed under $\ell_{p'}$-direct sum. Thus, $(U_\gamma^\pi,\ell^{p'}(G/H,E)) \in \text{Rep}_{p'}(G)$ and hence, $$(\text{Ind}_{H \uparrow G}(\pi),L^{p'}(G,H,\pi,E)) \in \text{Rep}_{p'}(G).$$
	\end{rems}
	
	\begin{rems}\label{4.1.2}
		It should be noted that $E$ sits inside $\ell_{p'}-\bigoplus\limits_{v \in G/H}E.$ Call $v_0 = H$ the standard base point of $G/H.$ For $\xi \in E,$ we define $f_\xi \in \ell_{p'}-\bigoplus\limits_{v \in G/H} E$ as $$f_\xi(v) = \begin{cases} 
			\xi & \text{if} \, \, v = v_0 \\
			0 & \text{otherwise}. 
		\end{cases}$$
		Then, by the above computations, $Mf_\xi \in L^{p'}(G,H,\pi,E)$ and is given by
		$$Mf_\xi(x) = Mf_\xi(\gamma(v)h) = \pi(h)^{-1} f_\xi(v) = \begin{cases} 
			\pi(h)^{-1}\xi & \text{if} \, \, v = v_0 \\
			0 & \text{otherwise}.
		\end{cases}$$
	\end{rems}
	
	\begin{rems} $($\cite[Corollary, Pg. 348]{cengiz}$)$
		The dual of $\ell^{p'}(G/H,E) = \ell_{p'}-\bigoplus\limits_{v \in G/H}E$ is given by $\ell_{p}-\bigoplus\limits_{v \in G/H}E^* = \ell^p(G/H,E^*)$ via the duality relation
		$$\langle f,g \rangle := \sum_{v \in G/H} \langle f(v),g(v) \rangle,$$
		where $f \in  \ell^{p'}(G/H,E)$ and $g \in \ell^p(G/H,E^*).$
	\end{rems}
	
	Now, consider the following Banach space:
	$$L^{p}(G,H,\pi^t,E^*) := \left\{ f:G \to E^* \, , \begin{array}{c} f(xh) = \pi(h)^t f(x) \, \forall \, x\in G, h \in H \\ \text{and} \sum\limits_{xH \in G/H} \|f(x)\|^{p} < \infty \end{array} \right\},$$
	with the norm $$\|f\|_p := \left( \sum\limits_{xH \in G/H} \|f(x)\|^{p} \right)^{1/{p}},$$
	where $\pi(\cdot)^t$ is the conjugate map from $E^*$ to $E^*.$
	
	By repeating the same arguments as above, it is easy to verify that the map $\tilde{M}:\ell^p(G/H,E^*) \to L^{p}(G,H,\pi^t,E^*)$ given by
	$$\tilde{M}f(x) = \tilde{M}f(\gamma(v)h) = \pi(h)^t f(v) \hspace{1cm} (f \in \ell^{p}(G/H,E^*)),$$ is an isometric isomorphism. Further, with $v_0 = H$ as the standard base point of $G/H,$ one can observe that $E^*$ sits inside $\ell_{p}-\bigoplus\limits_{v \in G/H}E^*,$ as in Remark \ref{4.1.2}.\\
	
	Here is the main result of this subsection.
	
	\begin{thms}
		Let $H$ be an open subgroup of $G.$ Given $v \in B_p(H),$ there exists $u \in B_p(G)$ such that $u|_H = v$ and $\|u\|_G =\|v\|_H.$ In other words, the restriction map $r: B_p(G) \to B_p(H)$ given by $r(u) = u|_H$ is a surjective isometry.
	\end{thms}
	
	\begin{proof}
		Let $v \in B_p(H).$ Then there exists $(\pi,E) \in \text{Rep}_{p'}(H)$ such that for every $h \in H,$
		$$v(h) = \langle \pi(h) \xi,\phi \rangle$$
		$\text{for some} \, \xi \in E \, \text{and} \, \phi \in E^*.$
		Now, define the function $u: G \to \mathbb{C}$ as
		$$u(g) := \begin{cases} 
			v(g) & \text{if} \, \, g \in H\\
			0 & \text{otherwise}. 
		\end{cases}$$
		In order to prove that $u \in B_p(G),$ we invoke the theory of $p$-induced representations. Consider the space $L^{p'}(G,H,\pi,E)$ as mentioned above. It follows from Remark \ref{QSP} that $$(\text{Ind}_{H \uparrow G}(\pi),L^{p'}(G,H,\pi,E)) \in \text{Rep}_{p'}(G).$$
		\textbf{Claim:} $u(g) = \langle (\text{Ind}_{H \uparrow G}(\pi))(g)\xi,\phi \rangle$ where $\xi \in E \hookrightarrow \ell_{p'}-\bigoplus\limits_{v \in G/H}E \cong L^{p'}(G,H,\pi,E)$ and $\phi \in E^* \hookrightarrow \ell_{p}-\bigoplus\limits_{v \in G/H}E^* \cong L^{p}(G,H,\pi^t,E^*).$
		
		Observe that
		\begin{align*}
			\langle (\text{Ind}_{H \uparrow G}(\pi))(g)\xi,\phi \rangle
			&= \langle (\text{Ind}_{H \uparrow G}(\pi))(g) Mf_\xi,\tilde{M}g_\phi \rangle\\
			&= \langle L_g(Mf_\xi),\tilde{M}g_\phi \rangle\\
			&= \sum\limits_{g'H \in G/H} \langle (L_g(Mf_\xi))(g'),\tilde{M}g_\phi(g') \rangle\\
			&= \sum\limits_{g'H \in G/H} \langle Mf_\xi(g^{-1}g'),\tilde{M}g_\phi(g') \rangle.\\
		\end{align*}
		Now, for a fix choice of a section $\gamma,$ the elements $g^{-1}g'$ and $g'$ can be uniquely written as
		$$g^{-1}g' = \gamma(q(g^{-1}g'))[\gamma(q(g^{-1}g'))^{-1} \, g^{-1}g']$$
		and 
		$$g' = \gamma(q(g'))[\gamma(q(g'))^{-1} \, g'].$$
		
		Observe that a non-zero term in the above summation appears if and only if $$q(g^{-1}g') = q(g') = v_0 (= H).$$ This forces us to conclude that $g$ and $g' \in H.$ Hence, to simplify the calculations, it is enough to take $g' = e$ in the above summation.
		
		Therefore, for $g \in H,$
		\begin{align*}
			\langle (\text{Ind}_{H \uparrow G}(\pi))(g)\xi,\phi \rangle &= \langle Mf_\xi(g^{-1}),\tilde{M}g_\phi(e) \rangle\\
			&= \langle \pi(g^{-1})^{-1} \xi,\phi \rangle\\
			&= \langle \pi(g) \xi,\phi \rangle = v(g) = u(g).
		\end{align*}
		For $g \notin H,$ it is clear that $\langle (\text{Ind}_{H \uparrow G}(\pi))(g)\xi,\phi \rangle = 0.$
		Thus, for each $v \in B_p(H),$ there exists $u \in B_p(G)$ such that $u|_H = v.$
		
		Now, it remains to show that the restriction map $r$ is an isometry. Clearly, for any $u \in B_p(G),$ we have
		$$\|u|_H\|_H \leq \|u\|_G.$$
		The non-trivial part is to establish the reverse inequality. Let $v \in B_p(H).$ Then, 
		$$v(g) = \sum_{n \in \mathbb{N}} \langle \pi_n(g)\xi_n,\phi_n \rangle,$$
		where $(\pi_n,E_n) \in \text{Cyc}_{p'}(H)$ with $\xi_n \in E_n, \phi_n \in {E_n}^*$ and $\sum\limits_{n \in \mathbb{N}} \|\xi_n\| \, \|\phi_n\| < \infty.$
		By repeating the above argument for each $\pi_n,$ we have the $p$-induced representation $\text{Ind}_{H \uparrow G}(\pi_n)$ for each $n \in \mathbb{N}.$ If we define $$u(g) := \sum\limits_{n \in \mathbb{N}} \langle (\text{Ind}_{H \uparrow G})(\pi_n)(g) Mf_{\xi_n},\tilde{M}g_{\phi_n} \rangle,$$ for $g \in G,$ then $u \in B_p(G)$ and $u|_H = v.$
		
		Now, by the definition of $u,$
		$$\|u\|_G \leq \sum\limits_{n \in \mathbb{N}} \|Mf_{\xi_n}\| \, \|\tilde{M}g_{\phi_n}\| = \sum\limits_{n \in \mathbb{N}} \|f_{\xi_n}\| \, \|g_{\phi_n}\| = \sum\limits_{n \in \mathbb{N}} \|\xi_n\| \, \|\phi_n\|.$$
		
		As a consequence of \cite[Lemma 4.3]{runde}, we have $\|u\|_G \leq \|v\|_H.$ 
	\end{proof}
	
	\subsection{For a compact normal subgroup $N$ of $G$}\label{sec42}
	Let $G$ be a locally compact group and $N$ be a normal compact subgroup of $G.$ The Haar measure on $N$ is normalized such that $|N| = 1.$ It is further assumed that the Haar measures $dg,d(gN)$ and $dn$ of $G,G/N$ and $N,$ respectively, are normalized such that the following formula holds for any $t \in C_c(G),$ i.e.,
	$$\int\limits_{G/N}\int\limits_{N} t(gn) dn \, d(gN) = \int\limits_{G} t(g) dg.$$
	
	\begin{nota}
		We write $h^g$ for $g^{-1}hg$ $(g,h \in G).$
	\end{nota}
	
	For any $u \in B_p(G),$ we define another function $u^g$ by
	$$u^g(h) = u(h^g) \hspace{1cm} (h \in G).$$
	
	Let us start with a basic lemma that relates the functions $u$ and $u^g,$ for $g \in G.$
	
	\begin{lems}\label{iso}
		For each $g \in G,$ the mapping $u \mapsto u^g$ from $B_p(G)$ to $B_p(G)$ is an isometry.
	\end{lems}
	
	\begin{proof}
		Let $u \in B_p(G)$ and $g \in G.$
		We need to prove that the function $u^g \in B_p(G)$ and $\|u^g\| = \|u\|.$
		Since $u \in B_p(G),$
		$$u(h) = \sum_{n \in \mathbb{N}} \langle \pi_n(h) \xi_n,\phi_n \rangle,$$
		where $(\pi_n,E_n) \in \text{Cyc}_{p'}(G)$ with $\xi_n \in E_n, \phi_n \in {E_n}^*$ and $\sum\limits_{n \in \mathbb{N}} \|\xi_n\| \, \|\phi_n\| < \infty.$
		Now, $$u^g(h) = u(g^{-1}hg) = \sum_{n \in \mathbb{N}} \langle \pi_n(g^{-1}hg) \xi_n,\phi_n \rangle = \sum\limits_{n \in \mathbb{N}} \langle \pi_n(h) (\pi_n(g)\xi_n),(\pi_n(g^{-1})^t\phi_n) \rangle.$$
		
		For each $n \in \mathbb{N},$ take ${\xi_n}' = (\pi_n(g)\xi_n) \in E_n$ and ${\phi_n}' = (\pi_n(g^{-1})^t\phi_n) \in {E_n}^*.$ Then,
		$$u^g(h) = \sum_{n \in \mathbb{N}} \langle \pi_n(h) {\xi_n}',{\phi_n}' \rangle \in B_p(G),$$
		and 
		$$\|u^g\| \leq \sum_{n \in \mathbb{N}} \|{\xi_n}'\| \, \|{\phi_n}'\| = \sum_{n \in \mathbb{N}} \|\pi_n(g) \xi_n\| \, \|\pi_n(g^{-1})^t \phi_n\| = \sum_{n \in \mathbb{N}} \|\xi_n\| \, \|\phi_n\|.$$
		
		Now, from \cite[Lemma 4.3]{runde}, it follows that $\|u^g\| \leq \|u\|.$ The reverse inequality is a consequence of the fact that $u = v^{g^{-1}}$ where $v = u^g.$
	\end{proof}
	
	The next lemma gives that the left and right translations act continuously and isometrically on $B_p(G).$ More precisely,
	
	\begin{lems}\label{lem422}
		For each $u \in B_p(G),$ the map $g \mapsto L_g(u)$ $(g \mapsto R_g(u))$ from $G$ to $B_p(G)$ is continuous and for each $g \in G,$ the map $u \mapsto L_g(u)$ $(u \mapsto R_g(u))$ is an isometry from $B_p(G)$ to $B_p(G),$ where $L_g(u)$ is left $(R_g(u)$ is right$)$ translation of $u.$
	\end{lems}
	
	\begin{proof}
		Let $u \in B_p(G)$ and $g \in G.$
		Then, $$u(h) = \sum_{n \in \mathbb{N}} \langle \pi_n(h) \xi_n,\phi_n \rangle$$
		where $(\pi_n,E_n) \in \text{Cyc}_{p'}(G)$ with $\xi_n \in E_n, \phi_n \in {E_n}^*$ and $\sum\limits_{n \in \mathbb{N}} \|\xi_n\| \, \|\phi_n\| < \infty.$
		
		It is enough to prove the lemma for left translation since the proof for right translation is similar. Observe that for each $g \in G,$ the map $u \mapsto L_g(u)$ is an isometry follows on the same lines as proved in Lemma \ref{iso}. Now, to prove that for each $u \in B_p(G),$ the map $g \mapsto L_g(u)$ is continuous, consider for $h \in G,$
		\begin{align*}
			(L_g(u)-u)(h) &= (L_g(u))(h)-u(h)\\
			&= \sum_{n \in \mathbb{N}} \langle \pi_n(g^{-1}h)\xi_n,\phi_n \rangle - \sum_{n \in \mathbb{N}} \langle \pi_n(h)\xi_n,\phi_n \rangle\\
			&= \sum_{n \in \mathbb{N}} \langle \pi_n(h)\xi_n,\pi_n(g^{-1})^t \phi_n \rangle - \sum_{n \in \mathbb{N}} \langle \pi_n(h)\xi_n,\phi_n \rangle\\
			&= \sum_{n \in \mathbb{N}} \langle \pi_n(h)\xi_n,(\pi_n(g^{-1})^t\phi_n-\phi_n )\rangle.\\
		\end{align*}
		Thus, 
		\begin{align*}
			\|L_g(u)-u\| &\leq \sum\limits_{n \in \mathbb{N}} \|\xi_n\| \, \|(\pi_n(g^{-1})^t\phi_n-\phi_n )\| \\ &\to 0 \, \text{as} \, g \to e.
		\end{align*}
		
		The last implication follows from the strong continuity of $\pi^t$ and the fact that the series $\sum\limits_{n \in \mathbb{N}} \|\xi_n\| \, \|\phi_n\|$ is convergent. Hence the result follows.
	\end{proof}
	
	Before proceeding to the main theorem of this subsection, we recall the $p$-induced representations for any subgroup $H$ of $G$ from \cite{jaming}.
	If $\pi$ is a strongly continuous isometric representation of the subgroup $H$ of $G$ in a Banach space $X,$ then $(\text{Ind}_{H \uparrow G}(\pi),L^{p'}(G,H,\pi,X))$ is the induced representation of $G$ on $L^{p'}(G,H,\pi,X).$\\
	For $1 < p' < \infty,$ the space
	$$L^{p'}(G,H,\pi,X) :=\left\{f:G \to X,
	\begin{array}{c}
		\forall \, \phi \in X^*, g \mapsto \langle f(g),\phi \rangle \, \text{is measurable and}  \\
		\forall g\in G, h\in H, f(gh) = \pi(h)^{-1} f(g)
	\end{array}
	\right\},$$
	is a Banach space with the norm
	$$\|f\|_{p'} := \left[ \int\limits_{G/H} \|f(g)\|^{p'} \, d(gH) \right]^{1/{p'}}.$$
	
	As in Section \ref{sec41}, by using the conjugate map $\pi(\cdot)^t,$ one can define another Banach space
	$$L^{p}(G,H,\pi^t,X^*):=\left\{f:G \to X^*,
	\begin{array}{c}
		\forall \, \gamma \in X^{**}, g \mapsto \langle f(g),\gamma \rangle \, \text{is measurable and}  \\
		\forall g\in G, h\in H, f(gh) = \pi(h)^t f(g)
	\end{array}
	\right\},$$
	with the norm
	$$\|f\|_{p} := \left[ \int\limits_{G/H} \|f(g)\|^{p} \, d(gH) \right]^{1/{p}}.$$
	
	The next lemma gives a relation between $L^{p'}(G,H,\pi,X)$ and $L^{p}(G,H,\pi^t,X^*).$
	
	\begin{lems}\label{lem423}
		For any $f \in L^{p}(G,H,\pi^t,X^*),$ the map $T_f: L^{p'}(G,H,\pi,X) \to \mathbb{C}$ given by $$T_f(t) := \int\limits_{G/H} \langle t(g),f(g) \rangle \, d(gH)$$ belongs to the dual of $L^{p'}(G,H,\pi,X).$
	\end{lems}
	\begin{proof}
		The well definedness of the map $T_f$ follows from the fact that for all $g \in G$ and $h \in H,$ we have, $t(gh) = \pi(h)^{-1} t(g)$ and $f(gh) = \pi(h)^t f(g).$ The result is just a consequence of the H\"{o}lder's inequality.
	\end{proof}
	For any $f \in C_c(G,X),$ it is proved in \cite{jaming} that the function $M_{p'}f$ given by
	$$M_{p'}f(g) := \int\limits_{H} [\pi(h)]f(gh) \, d(gH)$$
	belongs to $C_c(G,H,\pi,X).$ Similarly, it is easy to verify that for any $f \in C_c(G,X^*),$ the function $M_pf$ given by
	$$M_{p}f(g) := \int\limits_{H} [\pi(h)^{-1}]^tf(gh) \, d(gH)$$
	belongs to $C_c(G,H,\pi^t,X^*).$
	
	The idea used in the proof of the following result is motivated from \cite{cow}.
	\begin{thmss}\label{thm424}
		Suppose that $N$ is a compact normal subgroup of $G.$ Then the restriction map $r: B_p(G) \to B_p(N)$ has dense range.
	\end{thmss}
	
	\begin{proof}
		It is trivial that $B_p(G)|_N \subseteq B_p(N).$ Thus, in order to establish that the range of the map $r$ is dense in $B_p(N),$ we prove that for any $v \in B_p(N)$ and $\epsilon > 0,$ there is some $u \in B_p(G)$ for which $\|u|_N- v\|_N < \epsilon.$\\
		Let $v \in B_p(N).$ Then there exists $(\pi,E) \in \text{Rep}_{p'}(N)$ such that for every $n \in N,$
		$$v(n) = \langle \pi(n) \xi,\phi \rangle,$$
		$\text{for some} \, \, \xi \in E \, \text{and} \, \phi \in E^*.$
		
		The crucial step is to find a function $u \in B_p(G).$ Given an arbitrary $\epsilon >0,$ using Lemma \ref{lem422}, one can find neighborhoods $U$ of $e$ in $G$ and $V$ of $e$ in $N$ such that 
		\begin{align*}
			\|v^g - v\|_N &< \epsilon/2 \, \, \forall \, g \in U \, \, \text{and}\\
			\|L_{n^{-1}}v - v\|_N &< \epsilon/2 \, \, \forall \, n \in V.
		\end{align*}
		Let $K$ be a compact neighborhood of $e$ in $G$ such that $K \subseteq U$ and $K^{-1}K \cap N \subseteq V.$ Take a non-negative continuous function $w$ on $G$ such that supp($w$) $\subseteq K$ and
		\begin{align*}
			1 &= \int\limits_{G/N}\left[\int\limits_{N} w(gn) dn\right]^2 d(gN).
		\end{align*} 
		
		Now, by repeating the argument as in \cite[Theorem 1]{cow} and using the functions $v$ and $w,$ the function $u$ on $G$ is defined as
		$$u(\tilde{g}) := \int\limits_{G} \int\limits_{N} w(gn) w(\tilde{g}g) v(n) dn \, dg \hspace{1cm} (\tilde{g} \in G).$$
		
		We prove that the function $u \in B_p(G).$ Observe that for any $\tilde{g}$ in $G,$
		\begin{align*}
			u(\tilde{g}) &= \int\limits_{G} \int\limits_{N} w(gn) w(\tilde{g}g) v(n) dn \, dg\\
			&= \int\limits_{G} \int\limits_{N} w({\tilde{g}}^{-1}gn) w(g) v(n) dn \, dg.
		\end{align*}
		Take $t(g) = \int\limits_{N} w(\tilde{g}^{-1}gn) w(g) v(n) dn.$
		Then, by using the following formula,
		$$\int\limits_{G/N}\int\limits_{N} t(g\tilde{n}) d\tilde{n} \, d(gN) = \int\limits_{G} t(g) dg,$$
		we have,
		\begin{align*}
			u(\tilde{g}) &= \int\limits_{G/N} \int\limits_{N} \int\limits_{N} w({\tilde{g}}^{-1}g\tilde{n}n) w(g\tilde{n}) \, \langle \pi(n)\xi,\phi \rangle \, dn \, d\tilde{n} \, d(gN)\\
			&= \int\limits_{G/N} \int\limits_{N} \int\limits_{N} w({\tilde{g}}^{-1}gn) w(g\tilde{n}) \langle \pi(\tilde{n}^{-1}n)\xi,\phi \rangle \, dn \, d\tilde{n} \, d(gN)\\
			&= \int\limits_{G/N} \int\limits_{N} \int\limits_{N} w({\tilde{g}}^{-1}gn) w(g\tilde{n}) \, \langle \pi(n)\xi,\pi({\tilde{n}}^{-1})^t\phi \rangle \, dn \, d\tilde{n} \, d(gN)\\
			&= \int\limits_{G/N} \left\langle \int\limits_{N} w({\tilde{g}}^{-1}gn) \pi(n)\xi \, dn, \int\limits_{N} w(g\tilde{n}) \pi({\tilde{n}}^{-1})^t\phi \, d\tilde{n} \right\rangle d(gN).
		\end{align*}
		
		Now, to prove that the function $u \in B_p(G),$ we invoke the theory of $p$-induced representations. Consider the induced representation Ind$_{N \uparrow G}(\pi),$ the representation of $G$ on $L^{p'}(G,N,\pi,E)$ defined by 
		$$(\text{Ind}_{N \uparrow G}(\pi))(\tilde{g})f(g) = (L_{\tilde{g}}(f))(g) = f(\tilde{g}^{-1}g).$$
		
		By \cite[Theorem 3.1]{jaming}, $L^{p'}(G,N,\pi,E)$ can be identified with a quotient of $L^{p'}(G,E).$ Since $L^{p'}(G,E)$ is a QSL$_{p'}$-space by \cite[Theorem 3.1]{runde}, it follows that $L^{p'}(G,N,\pi,E)$ is a QSL$_{p'}$-space. As a consequence,
		$$(\text{Ind}_{N \uparrow G}(\pi),L^{p'}(G,N,\pi,E)) \in \text{Rep}_{p'}(G).$$
		
		Also, observe that for $\xi \in E$ and $\phi \in E^*,$ the functions given by $w_\xi(g) := w(g) \xi$ and $w_\phi(g) := w(g) \phi$ belongs to $C_c(G,E)$ and $C_c(G,E^*),$ respectively. Then, the functions $$M_{p'}w_\xi \in C_c(G,N,\pi,E) \subseteq L_{p'}(G,N,\pi,E), \text{where} \, M_{p'}w_\xi(g) = \int\limits_{N} w(gn) \pi(n)\xi \, dn ,$$ and $$M_pw_\phi \in C_c(G,N,\pi^t,E^*) \subseteq L_p(G,N,\pi^t,E^*), \text{where} \, M_pw_\phi(g) = \int\limits_{N} w(g\tilde{n}) \pi({\tilde{n}}^{-1})^t\phi \, d\tilde{n},$$ respectively. Thus, the function $u$ is given by
		$$u(\tilde{g}) = \left\langle \text{Ind}_{N \uparrow G}(\tilde{g})(M_{p'}w_\xi),M_pw_\phi \right\rangle.$$
		
		Hence, from Lemma \ref{lem423}, we have, $u \in B_p(G).$
		
		Now, it is left to prove that for this function $u \in B_p(G),$ we have $\|u|_N - v\|_N < \epsilon.$ This follows verbatim to the proof in \cite[Theorem 1]{cow} and hence, we omit it.
	\end{proof}
	
	\subsection{For a compact subgroup $H$ of an [SIN]$_H$-group G}\label{sec43}
	Let $G$ be a locally compact group with unit $e.$ Then $G$ is said to be an [SIN]$_H$-group if there exists a neighborhood base $\mathcal{V}$ at $e$ such that $h^{-1}Vh = V$ for all $V \in \mathcal{V}$ and $h \in H.$
	
	For $H = G,$ i.e., [SIN]$_G$-groups are simply referred to as SIN-groups. It is well known that [SIN]$_H$-groups are unimodular.
	
	Here is the main theorem of this subsection.
	
	\begin{thms}
		Let $H$ be a compact subgroup of an [SIN]$_H$-group $G.$ Then the range of the restriction map $r: B_p(G) \to B_p(H)$ is dense in $B_p(H).$
	\end{thms}
	\begin{proof}
		The main idea in proving this theorem is similar to the one used in Theorem \ref{thm424}. The proof follows exactly on the same lines as that of \cite[Theorem 2]{cow} and the final expression obtained for the function $u$ is
		$$u(\tilde{g}) = \int\limits_{G/H} \left\langle \int\limits_{H} w({\tilde{g}}^{-1}gh) \pi(h)\xi \, dh, \int\limits_{H} w(g\tilde{h}) \pi({\tilde{h}}^{-1})^t\phi \, d\tilde{h} \right\rangle d(gH),$$
		which is exactly the same as in the corresponding part of the proof in our Theorem \ref{thm424}. Thus, by repeating the argument as above, we have, $u \in B_p(G)$ which satisfies $\|u|_H - v\|_H < \epsilon$ for a given $v \in B_p(H)$ and $\epsilon > 0.$ This proves that the range of the map $r$ is dense in $B_p(H).$
	\end{proof}
	
	\section*{Acknowledgement}
	The first author is grateful to Indian Institute of Technology Delhi for the Institute Assistantship. The authors express their gratitude to Dr. Bharat Talwar for his valuable time and insightful suggestions that improved the manuscript.
	
	\section*{Data Availability}
	Data sharing does not apply to this article as no datasets were generated or analyzed during the current study.
	
	\section*{competing interests}
	The authors declare that they have no competing interests.

\end{document}